\documentstyle[12pt]{article}

  \newcommand{\const}{\rm const}
  \newcommand{\Var}{\rm Var}
  
  \newcommand{\Law}{\rm Law}

  \newcommand{\vraisup}{\rm vraisup}

  \newcommand{\MDT}{\rm  MDT}

\begin{document}

   \begin{center}

 {\bf  Non - asymptotic tails estimations for sums of random vectors} \par

\vspace{4mm}

{\bf having moderate decreasing tails.} \par

\vspace{5mm}

{\bf M.R.Formica,  E.Ostrovsky and L.Sirota}.\par

\end{center}

 Universit\`{a} degli Studi di Napoli Parthenope, via Generale Parisi 13, Palazzo Pacanowsky, 80132,
Napoli, Italy. \\

e-mail: mara.formica@uniparthenope.it \\

\vspace{4mm}

Department of Mathematics and Statistics, Bar-Ilan University, \\
59200, Ramat Gan, Israel. \\

e-mail: eugostrovsky@list.ru\\

\vspace{4mm}

\ Department of Mathematics and Statistics, Bar-Ilan University, \\
59200, Ramat Gan, Israel. \\

e-mail: sirota3@bezeqint.net \\

\vspace{5mm}

\begin{center}

{\sc Abstract.} \\

\end{center}

\vspace{4mm}

 \ We derive the sharp  non - asymptotical uniform  estimations for tails of distributions for classical normed sums
of centered normed independent random vectors having a moderate decreasing individual tails of summands. \par

\vspace{5mm}

  \par

\vspace{5mm}

\section{Statement of problem.}

\vspace{5mm}

 \hspace{3mm} Let $ \  B = (B, ||\cdot||B) \ $ be separable Banach space equipped with the norm $ \ ||\cdot||B; \ $
we will concentrate our  attention on the space of continuous numerical valued functions defined on certain
compact metric space. We will write  for the case when $ \ B \ $ is ordinary real line $ \ ||x||B = |x| \ - \ $ the ordinary
absolute value. \par
 \ Let $ \ (\Omega, \cal{B}, {\bf P} )  \ $ with expectation $ \ {\bf E} \ $ and variance $ \ \Var \ $
 be certain probability space and let
  also $ \ \xi \ $ be a centered in the weak sense  (mean zero) random variable (vector)  (r.v.)

 $$
  \forall y \in B^* \ \Rightarrow  {\bf E} y(\xi) = 0,
 $$
 with values in the space $ \ B: \ {\bf P} (\xi \in B) = 1 \ $ and define by $ \ \xi_i \ $  the independent copies of $ \ \xi. \ $ \par

 \ Let us  assume that the r.v. $ \ \xi \ $ has a finite weak second moment:

 $$
  \forall y \in B^* \ \Rightarrow  {\bf \Var} \ [ y(\xi)] < \infty.
 $$

 \ Put  therefore the normed as ordinary sum

\begin{equation} \label{Sn}
S_n := n^{-1/2} \sum_{i=1}^n \xi_i, \ n = 1,2,3,\ldots.
\end{equation}

 \ Further, define for arbitrary such a r.v. $ \ \xi \ $ its tail function

\begin{equation} \label{tail fun}
T[\xi](u) \stackrel{def}{=}  {\bf P} ( ||\xi||B > u), \ u > 0.
\end{equation}

 \vspace{3mm}

  \ The investigation of tail of distribution of the  sequence of the r.v. - s $ \  \{ S_n \} \ $ is the classical theme for
  the probability theory. We devote this preprint to the investigation of the  case when the r.v.- s $ \ \{\xi_i\} \ $ have a
  so - called  {\it moderate tails of distribution.} \par

 \vspace{4mm}

\ {\bf Definition 1.1.} We will say that the r.v. $ \ \xi \ $ has a moderate decreasing tail of distribution,
write $ \ \Law  (\xi) \in MDT, $ or equally write more detail

$$
\xi \in MDT(\beta, \gamma, V),
$$
if  for some constants $ \ \beta = \const > 2, \ \gamma = \const \in R,  \ $ and for some positive
continuous {\it slowly varying at infinity} function $ \ V = V(y), \ y \in (1,\infty) \ $

\begin{equation} \label{MDT}
T[\xi](u) \le u^{-\beta} \ \ln^{\gamma} u \ V(\ln u), \ u \ge e.
\end{equation}

 \vspace{4mm}

 \ {\bf  Our claim in this preprint is an exact uniform tail estimation for natural normed sums of random vectors having a moderate
 decreasing tails of distributions: }\\

\vspace{3mm}

\begin{equation} \label{Uniform tail}
Q(u) = Q[\xi](u) \stackrel{def}{=} \sup_n T[S_n](u), \ u \ge e,
\end{equation}

\vspace{3mm}

{\bf \ where } $ \ \xi \in MDT(\beta, \gamma, V). $ \par

\vspace{4mm}

 \ Notice, that it follows from the condition   (\ref{MDT}), in particular, $ \ \beta > 2, \ $ that $ \ {\bf \Var} (\xi)  < \infty; \ $
therefore the norming sequence $ \ 1/\sqrt{n} \ $  in (\ref{Sn}) \ is natural. \par

\vspace{3mm}

 \ Of course, these estimations may be used as ordinary in the statistics and in the method Monte - Carlo, see e.g.
\cite{Frolov Tchentzov}, \cite{Grigorjeva Ostrovsky}. In detail, let $ \ a \ $ be some unknown vector in the space $ \ B \ $
such that for certain r.v. $ \ \zeta \ $ having a distribution $ \ \mu: \ {\bf P}(\zeta \in A) = \mu(A) \ $ for any Borelian set $ \ A \ $
from the whole space $ \ B \ $  there holds

$$
a = {\bf E} g(\zeta) = \int_{B} g(w) \ \mu(dw).
$$
 \ Let $ \ \zeta_i, \ i = 1,2,\ldots \ $ be independent copies $ \ \zeta. \ $ The consistent a.e. as $ \ n \to \infty \ $  in the norm $ \ ||\cdot||B \ $
 estimate of the  value $ \ a \ $ has a form

$$
a_n := n^{-1} \sum_{i=1}^n g(\zeta_i).
$$
 \ Here $ \ \xi_i = g(\zeta_i) - a, \ $ and the estimation of the variable $ \ Q[g(\zeta)](u) \ $ as $ \ u \to \infty \ $ may be used for the building
 of non - asymptotical confidence region for the value $ \ a \ $ in the norm $ \ ||\cdot||B. \ $ \par

 \vspace{4mm}

 \ Roughly speaking,  we will ground that if the source random variable $ \ \xi \ $ has the moderate decreasing tail of distribution,
 then  under appropriate conditions   the r.v. - s $ \ S_n $  have also the moderate decreasing tail of distribution uniformly
 relative the parameter $ \ n. \ $ \par

\vspace{5mm}

\section{One dimensional case.}

\vspace{5mm}

\begin{center}

 \ {\sc  Some facts from the theory of Grand Lebesgue Spaces (GLS). } \par

\end{center}

\vspace{4mm}

 \ Define as ordinary for arbitrary numerical valued r.v. $ \ \zeta  \ $ its  Lebesgue - Riesz  $ \ L_p \ $ norm

$$
|| \zeta||_p \stackrel{def}{=} \left[ \ {\bf E} |\zeta|^p  \ \right]^{1/p}, \ p \in [1,\infty),
$$

 $$
 ||\zeta||_{\infty} \stackrel{def}{=} \vraisup_{\omega \in \Omega} |\eta(\omega)|.
 $$

 \vspace{3mm}

 \hspace{3mm}  Recall, see e.g.  \cite{Ahmed Fiorenza Formica at all}, \cite{Buldygin},
 \cite{fiokarazanalanwen2004}, \cite{fioguptajainstudiamath2008},
\cite{Fiorenza-Formica-Gogatishvili-DEA2018}, \cite{fioforgogakoparakoNAtoappear}, \cite{fioformicarakodie2017},
\cite{Kozachenko1},   \cite{Liflyand}, \cite{Naimark Ostrovsky}, \cite{Osekovski}, \cite{Ostrovsky1}, \cite{Ostrovsky 2004},
\cite{Ostrovsky HIAT} and so one,  that the so - called Grand Lebesgue Space
  (GLS)  $ \  G\psi = G\psi(b), \ b = \const \in (2,\infty] \ $  equipped with the norm  $ \ ||\zeta|| G\psi \ $  of the r.v. $ \ \zeta \ $ is defined as follows:
$ \ G\psi = \{\zeta: \ \Omega \to R, \ ||\zeta||G\psi < \infty, \} \ $ where

$$
||\zeta||G\psi \stackrel{def}{=} \sup_{p \in[2,b)}  \left[ \ \frac{||\zeta||_p}{\psi(p)} \ \right].
$$
 \ Here $ \ \psi = \psi(p), \ p \in [2,b)  \ $ is bounded from below measurable  function,
 which is names as ordinary as a
 {\it generating function} for this space. The set of all such a functions will be denoted by $ \ \Psi \ $ or more concrete $ \Psi(b). \ $ \par
 \ Note that one can consider a more general case $ \ p \in (1, \infty], \ $ but in this report only $ \ p \in [2,b), \ b \in (2,\infty). \ $ \par

\vspace{3mm}

 \ A very popular class of these spaces form the  {\it subgaussian} random variables, i.e. for
which $ \ \psi(p) = \psi_2(p) = \sqrt{p}, \ b = \infty. $ \par

 \ More generally,

$$
\psi(p) = \psi_m(p) := p^{1/m}, \ p \ge 1.
$$

 \ Suppose the r.v. $ \  \zeta  \ $  belongs to the space $ \ G \psi_m. \ $
The correspondent tail estimate is follow:

$$
\max \left[  {\bf P}(\zeta \ge u), \ {\bf P} (\zeta \le - u)  \right] \le \exp \left\{ -  (u/K)^m   \right\}, \ u > 0,
$$
and correspondent inverse conclusion also holds true.\par

 \ These space are used   in particular for obtaining of the exponential decreasing tail estimates for sums of
 random variables, independent or not, see e.g.
\cite{Kozachenko1},  \cite{Ostrovsky1},  sections 1.6, 2.1 - 2.5.  \par
 \ For instance, if $ \ {\bf E} \xi = 0 \ $  and  for some value $  \ m = \const > 0 $ and $ \ b = \infty \ $

$$
\max \left[  {\bf P}(\xi \ge u), \ {\bf P} (\xi \le - u)  \right] \le \exp \left\{ -  u^m   \right\}, \ u > 0,
$$
then

$$
 \sup_n \max \left[  {\bf P}(S(n) \ge u), \ {\bf P} (S(n) \le - u)  \right] \le
$$

$$
\exp \left\{ - C(m)  u^{\min(m,2)} \right\}, \ u > 0, \ C(m) \in (0, \infty),
$$
and the last estimate is essentially non - improvable. \par

 \vspace{3mm}

  \ Another possibility; the so - called {\it natural function.} Let $ \ \eta \ $ be a random variable (r.v.) such that

$$
 \exists b \in (2, \infty],  \ \forall p \in [2,b) \ \Rightarrow \ ||\xi||_p < \infty.
$$
 \ By definition, the natural function   $ \ \psi[\xi](p)  \  $ for the r.v. $  \ \xi \ $  is defined as follows

\begin{equation} \label{natural fun}
\psi[\xi](p)  \stackrel{def}{=} ||\xi||_p, \ p \in [2,b).
\end{equation}

 \ Evidently, $ \  ||\xi||G\psi[\xi] = 1. \ $ \par

\vspace{3mm}

 \ The belonging of some r.v. $ \ \xi \ $ to certain GLS $ \ G\psi, \ \psi \in \Psi \ $ is closely related with its tail behavior.
 Indeed, denote  $ \ \nu(p) := p \ln \psi(p), \ $

$$
\nu^*(y) := \sup_{p \in [2,b)}(py - \nu(p)), \ y \ge e, -
$$
the so - called (regional)  Young - Fenchel transform.  If the r.v. $ \ 0 \ne \xi \in G\psi, \ $ for definiteness let $ \ ||\xi||G\psi = 1, \ $
then

\begin{equation} \label{mom tail}
T[\xi](z) \le \exp \left( \ - \nu^*(\ln z)  \ \right), \ z \ge e.
\end{equation}

\vspace{3mm}

 \ Conversely, let the tail function $ \ T[\xi](x) \ $ for the one - dimensional random variable $ \ \xi \ $ be given; then

$$
{\bf E} |\xi|^p = p \int_0^{\infty} x^{p-1} \ T[\xi](x) \ dx,
$$
and we observe

\begin{equation} \label{norm from tail}
||\xi||G\psi = \sup_p \left\{ \ \frac{\left[ \ p \int_0^{\infty} x^{p-1} \ T[\xi](x) \ dx \ \right]^{1/p} }{\psi(p)}  \ \right\}.
\end{equation}

 \vspace{4mm}

 \begin{center}

 {\sc Examples.}

 \end{center}

 \vspace{4mm}

  \hspace{3mm}    Assume that the source one - dimensional r.v. $ \ \xi \ $ has a moderate decreasing tail
 of distribution,  see (\ref{MDT}). We intent to evaluate its natural function as $ \ p \in [2,\beta), \ $  which we will denote by
$ \  \kappa(p) = \kappa[\beta,\gamma,V](p). \ $ Namely,

$$
\kappa^p[\beta,\gamma,V](p) - e \stackrel{def}{=}  p^{-1} {\bf E} |\xi|^p - e \le \int_e^{\infty} x^{p - \beta -1} \  \ln^{\gamma}(x) \ V(\ln x) \ dx.
$$

\vspace{3mm}

 \ We need in this purpose to introduce the following auxiliary function  $ \ \theta[\gamma](p) =  \theta[\gamma]_{\beta, V}(p) \ $ for the values
 $ \ \ \gamma \in R, \ p \in [2, \beta) \ $

\vspace{3mm}

\ {\bf A.}  \ $ \ \gamma > - 1 \ \Rightarrow \ $

\begin{equation} \label{theta1}
\theta[\gamma](p) =  \theta[\gamma]_{\beta, V}(p) \stackrel{def}{=}  (\beta - p)^{-\gamma - 1} \ V(1/(\beta - p)).
\end{equation}

\vspace{3mm}

\ {\bf B.}  \ $ \ \gamma = - 1 \ \Rightarrow \ $

\begin{equation} \label{theta2}
\theta[-1](p) = \theta[\gamma](p) =  \theta[\gamma]_{\beta, V}(p) \stackrel{def}{=} |\ln (\beta - p)| \ V(1/(\beta - p)).
\end{equation}

\vspace{3mm}

\ {\bf C.}  \ $ \ \gamma < - 1 \ \Rightarrow \ $

\begin{equation} \label{theta2}
\theta[\gamma](p) =  \theta[\gamma]_{\beta, V}(p) \stackrel{def}{=} V(1/(\beta - p)).
\end{equation}

\vspace{4mm}

 \ {\bf Proposition 2.1.}

 \vspace{3mm}

\begin{equation} \label{comparison}
\kappa^p(p) = \kappa^p[\beta,\gamma,V](p) \asymp \theta[\gamma]_{\beta, V}(p), \ p \in [2,\beta).
\end{equation}

\vspace{4mm}

 \ {\bf Proof.}  \ We must  consider separately the three cases. \par

 \vspace{3mm}

 \hspace{3mm} {\bf A.} $ \ \gamma > - 1. \ $  \par

 \vspace{3mm}

  \ Let the r.v. $ \ \xi \ $ be as in (\ref{MDT}); we intend to estimate its (absolute) moment $  \ {\bf E} |\xi|^p, \ p \in [2,\beta). \ $
 We have as $ \ p \to \beta - 0 \ $

$$
 \kappa^p(p) \le  \int_1^{\infty} x^{p-1 - \beta} \ \ln^{\gamma} x \ V(\ln x) \ dx  =
$$

$$
\int_0^{\infty} e^{-y(\beta - p)} \ y^{\gamma} \ V(y) \ dy  = (\beta - p)^{-\gamma - 1} \ \int_0^{\infty} e^{-z} \ z^{\gamma} \ V(z/(\beta - p)) \ dz \sim
$$

$$
\Gamma(\gamma + 1) \ (\beta - p)^{-\gamma - 1} \ V(1/(\beta - p)) = \theta[\gamma](p).
$$
 where as ordinary $ \ \Gamma(\cdot) \ $ is Gamma function. See  for the more detail explanation \cite{Liflyand}, \cite{Ostrovsky HIAT}.\par

\vspace{3mm}

\hspace{3mm} {\bf B.} $ \ \gamma = - 1. \ $  Then   as $ \ p \in [2,\beta), \ p \to \beta - 0 \ $

\vspace{3mm}

$$
\kappa^p(p) \sim \int_e^{\infty} x^{p - 1 - \beta} \ \ln^{-1}(x) \  V(\ln x) \ dx =
\int_1^{\infty} e^{-y(\beta - p)} \ y^{-1} \ V(y) \ dy \sim
$$

$$
\int_{C (\beta - p)}^{\infty} e^{-z} \ z^{-1} \ V(z/(\beta - p)) \ dz \sim
$$

$$
V(1/(\beta - p)) \ \int_{C(\beta - p)}^{\infty} e^{-z} \ z^{-1} \ dz \sim V(1/(\beta - p)) \ \int_{(\beta - p)}^1 \ z^{-1} \ dz =
$$

$$
V(1/(\beta - p)) \ |\ln(\beta - p)| =\theta[-1](p).
$$

\vspace{3mm}

 \hspace{3mm} {\bf C.} $ \ \gamma < - 1. \ $  Then again  as $ \ p \in [2,\beta), \ p \to \beta - 0 \ $

\vspace{3mm}

$$
\kappa^p(p) \sim \int_e^{\infty} x^{p - 1 - \beta} \ \ln^{\gamma} x \ V(\ln x) \ dx =
\int_1^{\infty} e^{-(\beta - p)y} \ y^{\gamma} \ V(y) \ \ d y =
$$

$$
(\beta - p)^{-\gamma - 1} \int_{\beta - p}^{\infty} e^{-z} \ z^{\gamma} \  V(z/(\beta - p)) \ dz \sim
$$

$$
V(1/(\beta - p)) \ (\beta - p)^{-\gamma - 1}\int_{\beta - p}^e z^{\gamma} \ dz \asymp V(1/(\beta - p)) = \theta[\gamma](p).
$$

\vspace{5mm}

 \ {\bf Theorem 2.1.} \ Assume that the centered random variable satisfies the inequality (\ref{MDT}).
Denote $ \ \tau(p) := \ln \theta[\gamma]_{\beta,V}(p), \ 2 \le p < \beta.  \ $ Our proposition:

\vspace{4mm}

\begin{equation} \label{main d1}
Q[\xi](Cu) \le \exp \left\{ \  - \tau^*(\ln u) \  \right\}, \ u \ge e, \ C = C(\gamma,\beta,V) \in (0,\infty).
\end{equation}

\vspace{4mm}

 \ {\bf Proof.} We apply first of all  Proposition 2.1:

 \vspace{3mm}

\begin{equation} \label{moment estimation individual}
{\bf E} |\xi|^p \le C(\beta,\gamma,V) \cdot  \theta[\gamma]_{\beta, V}(p), \ p \in [2,\beta).
\end{equation}

\vspace{3mm}

 \ Further, we will use the famous moment estimations for the sums of the centered independent random variables, see e.g.
\cite{Dharmadhikari Jogdeo}, \cite{Ibragimov1},  \cite{Ibragimov2}, \cite{Naimark Ostrovsky}, \cite{Rosenthal}.
 Since the interval of the values $ \ p \ $ is bounded: $ \ 2 \le p < \beta, \ $ one can write

 \vspace{3mm}

\begin{equation} \label{moment estimation sums}
\sup_n {\bf E} |S_n|^p \le C_1(\beta,\gamma,V) \cdot  \theta[\gamma]_{\beta, V}(p), \ p \in [2,\beta).
\end{equation}

\vspace{3mm}

 \ It remains to apply the estimate  (\ref{mom tail}). \par

 \vspace{4mm}

 \begin{center}

 {\sc Examples.}

 \end{center}

 \vspace{4mm}

 \hspace{3mm}  {\bf I.}  Suppose  for example that the source one - dimensional centered r.v. $ \ \xi \ $ has a moderate decreasing tail
 of distribution, where  as above  $ \ \beta > 2, \ $ and suppose here that $ \gamma > -1, \ $  see (\ref{MDT}). Then

\vspace{3mm}

\begin{equation} \label{one dim result}
Q(u) = Q[\xi](u) \stackrel{def}{=} \sup_n T[S_n](u) \le C_1(\beta,\gamma,L) \  u^{-\beta} \ \ln^{\gamma + 1} u \ V(\ln u), \ u \ge e,
\end{equation}
and this estimation is  essentially  non - improvable.\par

\vspace{4mm}

 \ The non - improvability may be ground by means of consideration of the following example

\begin{equation} \label{counter MDT}
T[\xi](u) = u^{-\beta} \ \ln^{\gamma} u \ V(\ln u), \ u \ge e,
\end{equation}
as long as

$$
Q[\xi](u) \ge T[\xi](u) = u^{-\beta} \ \ln^{\gamma} u \ V(\ln u), \ u \ge e.
$$

\vspace{3mm}

 \ More complicated examples may be found in \cite{Liflyand}, \cite{Ostrovsky HIAT}. \par

\vspace{3mm}

 \ {\bf II.}  Suppose now that $ \ \gamma = - 1: \ $

$$
T[\xi](u) \le  u^{-\beta} \ \ln^{-1} u \ V(\ln u), \ u \ge e,
$$

then

$$
Q[\xi](u) \le C_2(\beta,V) \  \ln \ln u \ V(\ln u), \ u \ge e^e.
$$

\vspace{3mm}

\ {\bf III.} Let now $ \ \gamma <  - 1;   \ $ and suppose in addition  that
$ \  \lim_{u \to \infty} V(u) = 0. \ $ We state

$$
Q[\xi](u) \le C_3(\beta, \gamma, V) \ V(\ln u), \ u \ge e.
$$

\vspace{3mm}

 \ The possible lower bound in the both of  last examples is  quite  alike one in the first example:

$$
Q[\xi](u) \ge  T[\xi](u), \ u \ge e^e.
$$

\vspace{5mm}

\section{Main  result. Space of continuous functions.}

\vspace{5mm}

 \hspace{3mm} Let $ \  Z = \{z\}  \ $ be arbitrary set; the semi - distance function  $ \  \rho = \rho(z_1,z_2), \ z_{1,2} \in Z \ $ on this set will be
 clarified below. Recall that the semi - distance function is non - negative symmetrical function vanishing in the diagonal $ \ \rho(z,z) = 0, \ $
  satisfying the triangle inequality but  in general case the relation $ \ \rho(z_1,z_2) = 0 \ $ does not imply that  $ \ z_2 = z_1. \ $ \par

 \vspace{3mm}

 \ Let  $ \  \eta = \eta(z),  \ z \in Z  \ $ be separable centered: $ \ {\bf E} \eta(z) = 0 \ $ numerical valued random field (r.f.).
 Let also $ \ \eta_i = \eta_i(z), \ i = 1,2,\ldots \ $ be independent copies of $ \ \eta(z). \ $  Put as above

 $$
 Y_n(z) \stackrel{def}{=} n^{-1/2} \sum_{i=1}^n \eta_i(z)
 $$
and

$$
W(u) \stackrel{def}{=} {\bf P} (||Y_n(\cdot)|| > u), \ u \ge  e.
$$
 \ Hereafter the symbol $ \ ||f(\cdot)|| \ $ denotes an uniform norm of the function $ \ f: \ $

$$
||f|| = ||f(\cdot)|| \stackrel{def}{=} \sup_{z \in Z} |f(z)|.
$$

 \ The space of all the {\it continuous} numerical valued functions $ \ f: Z \to R \ $ equipped with the
 uniform norm $ \ ||\cdot|| \ $  will be denoted  as usually by $ \ C(Z) = C(Z,\rho). \ $ \par

\vspace{3mm}

 \ Let us suppose that

\begin{equation} \label{Uniform MDT}
\sup_{z \in Z} T[\eta(z)](u) \le u^{-\beta} \ \ln^{\gamma} u \ V(\ln u), \ u \ge e \ -
\end{equation}
the uniform \ \MDT \  condition.\par

  \ Define the following  {\it natural}  generating function

\begin{equation} \label{uniform psi}
\psi_{\gamma}(p) =  \psi[\gamma]_{\beta,V}(p) \stackrel{def}{=} \theta^{1/p} [\gamma]_{\beta,V}(p), \ 2 \le p < \beta,
\end{equation}
then it follows from (\ref{uniform psi}) that

\begin{equation} \label{boundedness psi norm}
\sup_{z \in Z} ||\eta(z)||G\psi[\gamma]_{\beta,V} = C_5 < \infty.
\end{equation}

 \ Let us introduce then the following bounded {\it natural} distance, more precisely, semi - distance on the set  $ \ Z \ $

\begin{equation} \label{natural distance}
 \rho(z_1,z_2) \stackrel{def}{=} ||\eta(z_1) - \eta(z_2)||G\psi[\gamma]_{\beta,V}, \ z_{1}, z_{2} \in Z.
\end{equation}

\vspace{4mm}

 \ Denote by $ \ H(\epsilon) = H(Z,\rho,\epsilon), \ 0 < \epsilon \le C_5, \ $ the metric entropy of the whole set (space) $ \ Z \ $
relative the metric $ \ \rho, \ $  i.e. the  (natural) logarithm of the minimal amount of the  closed ball in the distance $ \ \rho, \ $
which cover this set $ \ Z. \ $  Set $ \ N(\epsilon) = N(Z,\rho,\epsilon) = \exp H(Z,\rho,\epsilon). \ $ \par

\vspace{4mm}

 \ {\bf Theorem 3.1.} Suppose that $ \ \gamma > - 1 \ $ and  that the following {\it entropic integral} convergent:

 \vspace{3mm}

\begin{equation} \label{entropic int}
I(N) \stackrel{def}{=} \int_0^{C_5} N^{ (\gamma + 1)/\beta}(Z,\rho,\epsilon) \ d \epsilon < \infty.
\end{equation}

 \vspace{3mm}

  \ Then the set $ \ Z \ $  is a pre - compact semi - metric space relative to  the distance function $ \ \rho(\cdot, \cdot); \ $
  the random field $ \ \eta(z), \ $ as well as all the r.f. - s  $ \ Y_n(z) \ $  are $ \ \rho \ - \ $  continuous almost everywhere:

\vspace{3mm}

\begin{equation} \label{condit contin}
{\bf P} (\eta(\cdot) \in C(Z,\rho)) = {\bf P} (Y_n(\cdot) \in C(Z,\rho)) = 1
\end{equation}
and moreover

\begin{equation} \label{unif est}
\sup_n {\bf P} (\sup_{z \in Z} |Y_n(z)| > u) \le C_6(\gamma,\beta,V; I(N)) \ u^{-\beta} \ \ln^{\gamma + 1} u \ V(\ln u), \ u \ge e.
\end{equation}

\vspace{4mm}

 \ {\bf Proof.} It follows from the  proof of proposition of Theorem 2.1., indeed, we use the inequality (\ref{moment estimation sums}),
that uniformly relative both the parameters $ \ (z,n) \ $

\begin{equation} \label{moment estimation param}
\sup_n {\bf E} \ \sup_{z \in Z} |Y_n(z)|^p \le C_7(\beta,\gamma,V) \cdot  \theta[\gamma]_{\beta, V}(p) = C_7 \psi^p_{\gamma}(p), \ p \in [2,\beta),
\end{equation}

or equally

\begin{equation} \label{mom estimation with param}
\sup_n \ \sup_{z \in Z} ||Y_n(z)||_p \le C_8(\beta,\gamma,V) \cdot  \psi_{\gamma}(p), \ p \in [2,\beta);
\end{equation}

\begin{equation} \label{Gpsi estimation param}
\sup_n \ \sup_{z \in Z} ||Y_n(z)||G\psi_{\gamma} \le C_8(\beta,\gamma,V) < \infty.
\end{equation}

 \ We find quite analogously

\begin{equation} \label{Gpsi difference estimation}
\sup_n \ \sup_{z_1, z_2 \in Z, \ \rho(z_1,z_2) > 0 }   \left[ \  \frac{||Y_n(z_1) - Y_n(z_2)||G\psi_{\gamma}}{\rho(z_1, z_2)} \ \right] \le C_8(\beta,\gamma,V) < \infty.
\end{equation}

\vspace{3mm}

 \ Both the propositions of theorem 3.1 (\ref{condit contin}) and  (\ref{unif est}) follows immediately from Theorem 3.17.1
of monograph   \cite{Ostrovsky1}, chapter 3, section 17; see also an article \cite{Pisier}. \par

\vspace{4mm}

 \hspace{3mm} {\bf  Example 3.1.} Assume that $ \ Z = D \ $ is bounded convex  subset of the whole Euclidean space
 $ \ R^d \ $ equipped with ordinary Euclidean norm $ \ |z|, z \in D \subset R^d.\ $ Suppose that the distance $ \ \rho \ $
 is  such that

$$
\rho(z_1,z_2) \le C_9 \ |z_1 - z_2|^{\alpha},  \ \alpha = \const \in (0,1].
$$
 \ Then

$$
N(Z,\rho, \epsilon) \le C_{10} \ \epsilon^{-d/\alpha}, \ \epsilon \in (0,C).
$$
 \ The condition (\ref{entropic int}) is satisfied iff

$$
\frac{\beta}{\gamma + 1} > \frac{d}{\alpha}.
$$

\vspace{5mm}

\section{Concluding remarks.}

\vspace{5mm}

 \ Note that under conditions of theorem 3.1 the r.f. - s $ \ Y_n(\cdot) \ $ not only  are continuous a.e., but
satisfies the Central Limit Theorem in the space $ \ C(Z,\rho). \ $ This implies that the distributions of the r.f. $ \ Y_n(z) \ $
converges as $ \ n \to \infty \ $  weakly in this space to the distribution of the  centered Gaussian r.f. $ \ Y_{\infty}(z) \ $
having at the same covariation function as r.f. $ \ \eta(z). \ $  See \cite{Ostrovsky1}, chapter 4, section 4.4. \par

\vspace{6mm}

\vspace{0.5cm} \emph{Acknowledgement.} {\footnotesize The first
author has been partially supported by the Gruppo Nazionale per
l'Analisi Matematica, la Probabilit\`a e le loro Applicazioni
(GNAMPA) of the Istituto Nazionale di Alta Matematica (INdAM) and by
Universit\`a degli Studi di Napoli Parthenope through the project
\lq\lq sostegno alla Ricerca individuale\rq\rq .\par


\begin{thebibliography}{44}




\bibitem{Ahmed Fiorenza Formica at all}
{\bf I. Ahmed, A. Fiorenza, M.R. Formica, A. Gogatishvili, J.M. Rakotoson.}
{\it Some new results related to Lorentz G-Gamma spaces and interpolation.}
J. Math. Anal. Appl., {\bf 483,} \ {\bf 2}, (2020),  363  \ - \ 385.


\bibitem{Buldygin}
 {\bf V.~V.~Buldygin, D.~I.~Mushtary, E.~I.
~Ostrovsky} and {\bf M.~I.~Pushalsky.} {\it New Trends in
Probability Theory and Statistics.} Mokslas (1992), V.1, 78 \ - \ 92;
Amsterdam, Utrecht, New York, Tokyo.



\bibitem{Frolov Tchentzov}
{\bf Frolov A.S., Tchentzov N.N.}  {\it On the calculation by the Monte-Carlo
method definite integrals depending on the parameters.} Journal of Computational
Mathematics and Mathematical Physics, (1962), V. 2, Issue 4, p. 714 \ - \ 718, (in Russian.)



\bibitem{Grigorjeva Ostrovsky}
{\bf Grigorjeva M.L., Ostrovsky E.I.}  {\it Calculation of Integrals on discontinuous
Functions by means of depending trials method.} Journal of Computational
Mathematics and Mathematical Physics, (1996), V. 36, Issue 12, p. 28-39 (in Russian).





\bibitem{Dharmadhikari Jogdeo}
{\bf S.W.Dharmadhikari and K.Jogdeo.}  {\it Bounds on moments of certain random variables.}
 Ann. Math. Stat., {\bf 40,}  No. 4, 1506  \ - \ 1508, (1969).


\bibitem{fiokarazanalanwen2004}
{\bf A.~Fiorenza} and {\bf G.~E.~Karadzhov.} {\it Grand and small
Lebesgue spaces and their analogs}, Z. Anal. Anwendungen \textbf{23}
(2004), no.~4, 657--681.

\bibitem{fioguptajainstudiamath2008}
{\bf A.~Fiorenza, B.~Gupta} and {\bf P.~Jain.} {\it The maximal
theorem for weighted grand Lebesgue spaces}. Studia Math.
\textbf{188} (2008), no.~2, 123--133.


\bibitem{Fiorenza-Formica-Gogatishvili-DEA2018}
{\bf A.~Fiorenza, M.~R.~Formica} and {\bf A.~Gogatishvili.} {\it On
grand and small Lebesgue and Sobolev spaces and some applications to
PDE's}. \emph{Differ. Equ. Appl.} \textbf{10},  (2018), no.~1, 21--46.

\bibitem{fioforgogakoparakoNAtoappear}
{\bf A.~Fiorenza, M. R.~Formica, A.~Gogatishvili, T.~Kopaliani} and
{\bf J.~M. Rakotoson.} {\it Characterization of interpolation
between grand, small or classical Lebesgue spaces}. Preprint
arXiv:1709.05892, Nonlinear Anal., {to appear}.

\bibitem{fioformicarakodie2017}
{\bf A.~Fiorenza, M.~R.~Formica} and {\bf J.~M. Rakotoson.} {\it
Pointwise estimates for {$G\Gamma$}-functions and applications}.
Differential Integral Equations \textbf{30} (2017), no.~11-12,
809--824.

\bibitem{formicagiovamjom2015}
{\bf M.~R. Formica} and {\bf R.~Giova.} {\it Boyd indices in
generalized grand Lebesgue spaces and applications}. Mediterr. J.
Math. \textbf{12} (2015), no.~3, 987--995.


\bibitem{Ibragimov1}
{\bf Ibragimov R., Sharachmedov Sh.}  {\it On the exact constant in the Rosenthal
Inequality.} Theory Probab. Appl., 1997, V. 42. p. 294 \ - \ 302.

\bibitem{Ibragimov2}
{\bf Ibragimov R., Sharachmedov Sh.} {\it The exact constant in the Rosenthal Inequality for sums random variables with mean zero.}
 Probab. Theory Appl., 2001, V. 46, 1, p. 127  \ - \  132.


\bibitem{Kozachenko1}
{\bf Yu.~V.~Kozachenko} and {\bf E.~I.~Ostrovsky.} {\it The Banach
Spaces of random variables of sub-Gaussian type.} of Probab. and
Math. Stat., \textbf{32} (1985), (in Russian). Kiev, KSU, 43 \ -\ 57.

\bibitem{Liflyand}
{\bf E. Liflyand, E. Ostrovsky and L. Sirota.}  {\it Structural properties of bilateral grand Lebesque
spaces.} Turkish J. Math. 34 (2010), no. 2, 207{219.



\bibitem{Naimark Ostrovsky}
{\bf Naimark B.,  Ostrovsky E.}
{\it Exact Constants in the Rosenthal Moment Inequalities for Sums of
independent centered Random Variables.} \\
arXiv:math/0411614v1 [math.PR] 27 Nov 2004

\bibitem{Osekovski}
{\bf Osekowski  A.} {\it A  Note  on  Burkholder  \ - \ Rosenthal  Inequality.}
Bull. Polish Academy of Science, Math., {\bf 60,} (2012), 177 \  -\ 185.

\bibitem{Ostrovsky1}
{\bf E.~Ostrovsky.} {\it Exponential estimates for random fields and
its applications.} 1999, OINPE, Moscow - Obninsk.

\bibitem{Ostrovsky 2004}
{\bf E.~Ostrovsky.} {\it Bide-side exponential and  moment
inequalities for tails  of  distribution of  Polynomial
Martingales.} \\
 arXiv: math.PR/0406532 v.1 Jun. 2004.

\bibitem{Ostrovsky HIAT}
{\bf Ostrovsky E. and Sirota L.}  \it Moment Banach spaces: theory and applications.} HIAT
Journal of Science and Engineering, C, Volume 4, Issues 1 - 2, pp. 233 - 262, (2007).

\bibitem{Ostrovsky-Sirota Jan 2008}
{\bf E.~Ostrovsky} and {\bf L.~Sirota.} {\it  Exponential bounds in
the law of iterated logarithm for martingales.} \\
 arXiv:0801.2125v1 [math.PR], 14 Jan 2008.

\bibitem{Ostrovsky-Sirota-boundedeness operator bilateral GLS Oct2011}
{\bf E.~Ostrovsky} and {\bf L.~Sirota.} {\it Boundedness of
operators in bilateral Grand Lebesgue Spaces, with exact and weakly
exact constant calculation}. \\
 arXiv:1104.2963 [math.FA] Apr 2011.

\bibitem{Ostrovsky 2012}
{\bf E. Ostrovsky, L.Sirota.} {\it Moment and tail estimates for martingales and martingale transform,
with application to the martingale limit theorem in Banach spaces.} \\
arXiv:1206.4964v1 [math.PR] 21 Jun 2012

\bibitem{Ostrovsky-Sirota Oct2013}
{\bf E.~Ostrovsky} and {\bf L.~Sirota} {\it Simplification of the
majorizing measures method, with development.}  \\
arXiv:1302.3202v1 [math.PR]  13 Feb 2013.

\bibitem{Ostrovsky 2014}
{\bf E. Ostrovsky, L.Sirota.} {\it Sharp moment estimates for polynomial martingales.} \\
arXiv:1410.0739v1  [math.PR]  3 Oct 2014

\bibitem{Ostrovsky-Sirota Oct2015}
{\bf E.~Ostrovsky} and {\bf L.~Sirota.} {\it Vector rearrangement
invariant Banach spaces of random variables with exponential
decreasing tails of distributions.} \\
 arXiv:1510.04182v1 [math.PR] 14 Oct 2015.

\bibitem{Ostrovsky-Sirota-fundamental function GLS 2015}
{\bf E.~Ostrovsky} and {\bf L.~Sirota.} {\it Fundamental function
for Grand Lebesgue Spaces}. \\
 arXiv:1509.03644 [math.FA] Sept. 2015.

\bibitem{Pisier}
{\bf G. Pisier.}  {\it Conditions d'entropie assurant la continuitie de certains processus et applications
la l'analyse harmonique.} (French.) Seminaire d'analyse fonctionnelle, (1980), Exp. No. 13  \ - \ 14, pp. 43 \ - \ 46.

\bibitem{Rosenthal}
{\bf  Rosenthal H.P.}  {\it On the subspaces of} $ \ Lp \ (p >2) \ $ {\it spanned by sequences of independent variables.}
 Israel J. Math., 1970, \ No 3, p. 273 \  - \ 303.



\end{thebibliography}
\end{document}